\documentclass[11pt]{article}
\usepackage{amsthm,amsmath,amsfonts}
\usepackage[authoryear]{natbib}
\usepackage{color,graphicx,rotating}
\usepackage[mathscr]{euscript}
\usepackage{kbordermatrix}
\usepackage{setspace}
\usepackage{bm}
\def\wh#1{\widehat{#1}}
\def\bn{ \bm{n} }
\def\Qn{ Q_{ \bm{n} } }
\def\Qm{ Q_{ \bm{m} } }
\def\bx{\bm{x}}
\def\bX{\bm{X}}
\def\by{\bm{y}}

\def\bu{\bm{u}}
\def\bp{\bm{p}}
\def\bv{\bm{v}}
\def\bw{\bm{w}}
\def\bE{\mathbb{E}}
\def\C{ {\cal C} }
\def\G{ {\cal G } }
\def\KP {Krawtchouk polynomials }
\def\MKP {multivariate Krawtchouk polynomials }

\usepackage[sf,bf,tiny]{titlesec}
\titlelabel{\thetitle. \enspace}
\title{\normalsize\bf AN INTRODUCTION TO MULTIVARIATE KRAWTCHOUK POLYNOMIALS AND THEIR APPLICATIONS
}
\author{{\sc Persi Diaconis\thanks{Department of Statistics, Sequoia Hall, 
390 Serra Mall, Stanford University, Stanford, California 94305-4065, USA.}
 and Robert Griffiths\thanks{Department of Statistics, University of Oxford, 1 South Parks Rd, Oxford, OX1 3TG, UK; email  griff@stats.ox.ac.uk}
}\\\bigskip
\emph{Stanford University and Oxford University}
\\Version 1.2
 }

\begin{document}
\maketitle
\section*{Abstract}
Orthogonal polynomials for the multinomial distribution $m(\bm{x},\bm{p})$ of $N$ balls dropped into $d$ boxes (box $i$ has probability $p_i$) are called multivariate Krawtchouk polynomials. This paper gives an introduction to their properties, collections of natural Markov chains which they explicitly diagonalize and associated bivariate multinomial distributions.
\section*{Introduction}
This paper relates orthogonal polynomials, Markov chains and statistical modeling. Bill Studden loved all of these topics and their interactions.
The setting is a $d$-category multinomial distribution of $N$ balls dropped into $d$ boxes with probabilities 
$\bm{p}  = (p_1,\ldots ,p_d)$. Throughout $p_i > 0$ and the multinomial distribution is
\[
m(\bm{x},\bm{p}) = {N \choose x_1,\cdots ,x_d}\prod_{i=1}^dp_i^{x_i},\>\>\>0 \leq x_i,\>\sum_{i=1}^dx_i = N.
\]
\textcolor{black}{Systems of} orthogonal polynomials for the multinomial are defined in \citet{G1971}. They generalize Krawtchouk polynomials for the Binomial and are called multivariate Krawtchouk polynomials. A self contained development of these along with properties and pointers to an emerging literature is in section two. 

 \textcolor{black}{There is not a unique set of orthogonal polynomials in higher dimensions than one because the ordering of  
powers $x_1^{n_1}\cdots x_d^{n_d}$ is not unique within a fixed degree $|n| = n_1 + \cdots + n_d$. \textcolor{black}{The definitions here lead to constructing a unique set of polynomials by choosing an orthonormal basis of functions on $\{1,2,\ldots, d\}$, see the definitions in section 2.1.} Further orthogonal polynomials in a particular set are uniquely defined by their multiple leading coefficients} \textcolor{black}{(lemma 1 in section 2.1). This allows us to show that the polynomials suggested by \citet{GR2011} and \citet{X2013} fit into our definition.}

Multivariate Krawtchouk polynomials appear as the eigen-functions of a variety of natural Markov chains generalizing the classical Ehrenfest urn.
Consider an urn containing $N$ balls in $d$ colours. A ball is chosen at random and its colour changed to $j$ with probability $p_j$. This can be generalized in various directions. The balls can be partitioned into groups in a general way. Balls in the $\alpha^{\text{th}}$ group have their colours changed in a general way: a ball of colour $i$ is changed to colour $j$ with probability $P_{\alpha}(i,j)$. With appropriate choices, all of these chains are diagonalized by multivariate Krawtchouk polynomials. This work captures chains previously studied by \citet{HR2008}, \citet{KZ2009}, \citet{ZL2009} and    
\citet{M2010}, \citet{M2011}. Multivariate Krawtchouk polynomials also have a universal quality, diagonalizing symmetrized products of general Markov chains. These topics are explained in section three.

The third central topic is bivariate distributions with multinomial margins. This topic 
has a long history going back to work of \citet{L1969}.
Bivariate Lancaster distributions with multinomial margins
 have an expansion
 \[
 P(\bm{x},\bm{y}) =  m(\bm{x},\bm{p})m(\bm{y},\bm{p})
 \Big \{1 + \sum_{\bm{n}}\rho_{\bm{n}}h_{\bm{n}}Q_{\bm{n}}(\bm{x})Q_{\bm{n}}(\bm{y})\Big \},
 \]
 with $Q_{\bm{n}}$ the multivariate Krawtchouk polynomials, and $h_{\bm{n}}$ given by
 \[
 \mathbb{E}\big [Q_{\bm{n}}(\bm{X})Q_{\bm{n}}(\bm{Y})\big ]=\delta_{\bm{m}\bm{n}}h_{\bm{n}}\rho_{\bm{n}}.
 \]
 The $\rho_{\bm{n}}$ are called generalized correlations.
 A basic problem, the Lancaster problem, is what values of $\rho_{\bn}$ are admissible to have $P(\bm{x},\bm{y}) \geq 0$?.
 This problem was solved in the Binomial case by \citet{E1969}. A useful necessary and sufficient condition appears in section three. This leans on the multinomial hypergroup property which may be of independent interest. 
 There are natural choices of $\rho_{\bm{n}}$; if $K(\bm{x},\bm{y})$ is a reversible Markov chain with Krawtchouk polynomial eigenfunctions, then 
 \begin{equation*}
 P(\bm{x},\bm{y})=  m(\bm{x},\bm{p})K(\bm{x},\bm{y})
 \end{equation*}
 has a Lancaster expansion with $\rho_{\bm{n}}$ the eigen-values of $K$. This paper is a $d$-dimensional version of 
\citet{DG2012} which works out the connections between generalized Ehrenfest urns with two colours, Lancaster expansions and one variable Krawtchouk polynomials. The theory is more complete (and simpler) in this case and the reader might find it useful motivation.
 
 As usual, new developments raise new questions; what are the extreme points of the $d\times d$ stochastic matrices with $P$ as stationary distribution? Find a full solution of the Lancaster problem. What is the connection to Schur-Weyl duality and Bosonic Fock space? We are sorry not to be able to call on Bill Studden's expert help.
 \section*{2.1 Definitions and Background} The multinomial distribution associated with dropping $N$ balls into $d$ boxes having probabilities
${\bm p}  = (p_1,\ldots ,p_d)$ is
\begin{equation}
m(\bm{x},\bm{p}) = {N \choose x_1,\cdots ,x_d}\prod_{j=1}^dp_j^{x_j},\>\>\>0 \leq x_j,\>\sum_{j=1}^dx_j = N.
\label{s:2.1}
\end{equation}
\citet{G1971} defines orthogonal polynomials for $m(\bm{x},\bm{p})$ by choosing  a complete set of orthogonal functions $\{u^{(l)}_j\}$ with respect to $\bm{p}$, for $l = 0,1,\ldots,d-1$, $j = 1,2,\dots ,d$. Insist throughout that $u_j^{(0)} \equiv 1$. Thus, for all $k,l = 0,1,\ldots,d-1$,
\begin{equation}
\sum_{j=1}^du_j^{(k)}u_j^{(l)}p_j = \delta_{kl}a_k.
\label{s:2.2}
\end{equation}
\textcolor{black}{In this paper we usually take these functions to be orthonormal when $a_k=1$, $k=1,\ldots,d-1$, unless indicated.}
Examples of natural choices of $\{u^{(l)}\}$ are given in the following section 2.2. Often $\{u^{(l)}\}$ is a basis of eigen-functions for a Markov chain on $[d] = \{1,2,\ldots,d\}$. Writing $\bm{u}$ for $\{u^{(l)}\}$, \citet{G1971} defines a collection of orthogonal polynomials
\begin{equation}
\big \{Q_{\bm{n}}(\bm{x},\bm{u})\big \}\text{~with~}\bm{n} = (n_1,\ldots ,n_{d-1}),\>|\bm{n}| \leq N
\label{s:2.3}
\end{equation}
as the coefficient of $w_1^{n_1}\cdots w_{d-1}^{n_{d-1}}$ in the generating function
\begin{equation}
G(\pmb{x},\pmb{w}, \pmb{u}) 
= \prod_{j=1}^d\Big \{1 + \sum_{l=1}^{d-1}w_lu_{j}^{(l)}\Big \}^{x_j}.
\label{s:2.4}
\end{equation}
For integer $x_j$, expanding each term in the product by the binomial expansion gives a polynomial in $x_j$ so the coefficient of $\bm{w}^{\bm{n}}$ is a polynomial in $\bm{x}$. It is easy to see that $Q_{\bm n}$ is a polynomial of degree $|\bm {n}|= n_1+\cdots +n_{d-1}$.

\noindent
For example:
\begin{equation}
\begin{array} {lll}
\bm{n} &= (0,\ldots,0) &\Qn\equiv 1
\\
\bm{n} &= (0,\ldots,1_l,\ldots, 0) &\Qn= \sum_{j=1}^du_j^{(l)}x_j\overset{\text{Def}}= S_l
\\
\bm{n} &= (0,\ldots, 1_l,\ldots,1_m,\ldots ,0) &Q_{\bm{n}} 
=\frac{1}{2}S_lS_m - \frac{1}{2}\sum_{i=1}^dx_iu_i^{(l)}u_i^{(m)}
\\
\bm{n} &= (0,\ldots,2_l,\ldots, 0) &Q_{\bm{n}} = \frac{1}{2}S_l^2- \frac{1}{2}\sum_{i=1}^dx_i\big (u_i^{(l)}\big )^2.
\end{array}
\end{equation}
It is straightforward to show, using (\ref{s:2.4}), that the $\{\Qn\}$ are orthogonal:
\begin{equation}
\bE\big [\Qn(\bX,\bu)\Qm(\bX,\bu)\big ] 
= \delta_{\bm{m}\bm{n}}
{N\choose |\bm{n}| }{|\bm{n} | \choose \bm{n} }
\textcolor{black}{\prod_{j=1}^{d-1}a_j^{n_j}}.
\label{s:2.5}
\end{equation}
In (\ref{s:2.5}), $\bX$ has a multinomial $m(\bx,\bp)$ distribution.
\citet{G1971} also gives a related construction of multivariate orthogonal polynomials on the negative multinomial distribution via a generating function approach. 
 
The $\{\Qn\}$ also have an easily verified stability property: their definition does not depend on $N$, as long as $|\bn| \leq N$ the same
 $\{\Qn\}$ work for all sufficiently large $N$.
 
 As defined, it is not so clear how to express $\Qn$ as a polynomial. Recent work of \citet{MT2004} and \citet{GR2011} clarifies this. They use hypergeometric notation. Let
 \begin{equation}
F_1^{(n)}(-\bm{m},-\bm{x};-N;\bm{u}) :=
\sum_{k_{\cdot \cdot}\leq N}
\frac{
\prod_{i=1}^n(-m_i)_{(k_{i\cdot})}\prod_{j=1}^n(-x_j)_{(k_{\cdot j})}
}
{\prod_{ij}k_{ij}!(-N)_{(k_{\cdot\cdot})}} 
\prod_{i,j}u_{ij}^{k_{ij}},
\label{s:2.6}
\end{equation}
where a $\cdot$ in an index means sum, eg $k_{\cdot j} = \sum_{i=1}^{d-1}k_{ij}$,
\[
a_{(n)} = a(a+1)\cdots (a+n-1),
\>\>\>
a_{[n]} = a(a-1)\cdots (a-n+1)
\]
and the sum is over all $n\times n$ matrices $(k_{ij})$ with non-negative integer entries with sum of entries at most $N$.
\bigskip

\noindent
\textbf{Proposition 1.} Let $\{v_j^{(l)}\}$ be an orthonormal basis, 
$0 \leq l \leq d-1$, $1 \leq j \leq d$, with $v_d^{(l)}\equiv 1$.
For all $l$, let $\textcolor{black}{u_{ij} = 1 - v_j^{(i)}}$, $\textcolor{black}{i,j \in [d-1]}$, then \textcolor{black}{with} $\Qm(\bx,\bv)$ defined in (\ref{s:2.3}),(\ref{s:2.4}) 
\[
\Qm(\bx,\bv) = \frac{N!}{\prod_{i=1}^d m_i!} 
F_1^{(d-1)}(-\bm{m},-\bm{x};-N;\bm{u}).
\]
\bigskip

\noindent
\textbf{Proof} Reduce the variables to $x_1,\ldots,x_{d-1}$ by letting
$x_d = N - \sum_{i=1}^{d-1}x_i$, then
\begin{eqnarray*}
G(x,w) &=& \prod_{j=1}^d\Big (1 + \sum_{i=1}^{d-1}w_iv_{ j}^{(i)}\Big )^{x_j}
\nonumber \\
 &=& \prod_{j=1}^{d-1}\Big (1 + \sum_{i=1}^{d-1}w_iv_{ j}^{(i)}\Big )^{x_j}
\times \Big (1 + \textcolor{black}{w_{\cdot}}\Big )^{N-\sum_{j=1}^{d-1}x_j}
\nonumber \\
 &=& \prod_{j=1}^{d-1}\Big (1+ \textcolor{black}{w_{\cdot}} - \sum_{i=1}^{d-1}w_iu_{ij}\Big )^{x_j}
\times \Big (1 + \textcolor{black}{w_{\cdot}}\Big )^{N-\sum_{j=1}^{d-1}x_j}
\nonumber \\
 &=& \prod_{j=1}^{d-1}\Big (1 - \sum_{i=1}^{d-1}\frac{w_i}{1+ \textcolor{black}{w_{\cdot}}}u_{ij}\Big )^{x_j}
\times \Big (1 + \textcolor{black}{w_{\cdot}}\Big )^{N}
\nonumber \\
&=& 
\sum_{k_{ij}}\Big (1 + \textcolor{black}{w_{\cdot}}\Big )^{N-k_{\cdot\cdot}}\prod_{i=1}^{d-1}w_i^{k_{i\cdot}}
\prod_{ij}(-u_{ij})^{k_{ij}}
\prod_{j=1}^{d-1}\frac{{x_j}_{[k_{\cdot j}]}}{\prod_{i=1}^{d-1}k_{ij}!}.
\label{s:2.7}
\end{eqnarray*}
The coefficient of 
$
 \prod_{i=1}^{d-1}w_i^{m_i}
$
in (\ref{s:2.7}) is 
\begin{eqnarray*}
&&\sum_{k_{\cdot\cdot} \leq N}
\frac{{(N-k_{\cdot \cdot})!}}{\prod_{i=1}^{d-1}(m_i - k_{i\cdot})!}
\prod_{ij}(-u_{ij})^{k_{ij}}
\prod_{j=1}^{d-1}\frac{{x_j}_{[k_{\cdot j}]}}{\prod_{i=1}^{d-1}k_{ij}!}
\nonumber \\
&=&
\frac{N!}{\prod_{i=1}^{d-1}m_i!}
\sum_{k_{\cdot \cdot}\leq N}
\frac{
\prod_{i=1}^{d-1}(-m_i)_{(k_{i\cdot})}\prod_{j=1}^{d-1}(-x_j)_{(k_{\cdot j})}
}
{\prod_{ij}k_{ij}!(-N)_{(k_{\cdot\cdot})}} 
\prod_{i,j}u_{ij}^{k_{ij}}.
\end{eqnarray*}
Equating coefficients now gives the result.
\qed
\bigskip

A well known and easy to check conditional product binomial construction is
\begin{equation}
m(\bm{x};\bm{p}) = \prod_{j=1}^{d-1}b(x_j,p_j/(1-|\bm{p}_{j-1}|),N-|\bm{x}_{j-1}|),
\label{conditional:0}
\end{equation}
where $b(x,p,N) = {N\choose x}p^x(1-p)^{N-x}$, $\bm{x}_k = (x_1,\ldots x_k)$, $|\bm{x}_k|=x_1+\cdots +x_k$ with similar notation for the $\bm{p}_k$ terms.
The way that a set of multi-dimensional orthogonal polynomials is constructed using (\ref{conditional:0}) is to use 1-dimensional Krawtchouk orthogonal polynomials on the conditional distributions. 
\citet{X2013} gives a unified treatment of multivariate Hahn, Jacobi and Krawtchouk polynomials which are constructed from product conditional distributions. In his treatment there are natural constructions beginning with multivariate Jacobi polynomials which lead to Hahn polynomials via a multinomial-Dirichlet mixture and then from multivariate Hahn polynomials to multivariate Krawtchouk polynomials as a limit when the index parameters in the multinomial-Dirichlet tend to infinity with the ratios tending to $\bm{p}$. We show in the next theorem that his multivariate Krawtchouk polynomials are a special case of the multivariate Krawtchouk polynomials in this paper. A general preliminary lemma is needed.

\bigskip

\noindent
{\bf Lemma 1.} \emph{Let $\big \{Q_{\bm{n}}(\bm{x})\big \}$ be a $d$-dimensional orthogonal polynomials set on a random variable $\bm{X}$.
Then the polynomials are uniquely determined by their leading coefficients.
}
\bigskip

\noindent
{\bf Proof.} Take  the orthogonal polynomial set to be orthonormal without loss of generality. 
Denote the reproducing kernel polynomials by
\[
Q_n(\bm{x},\bm{y}) = \sum_{|\bm{n}| = n}Q_{\bm{n}}(\bm{x})Q_{\bm{n}}(\bm{y}),\>n=0,1,\ldots.
\]
The reproducing kernel orthogonal polynomials are invariant under all choices of multidimensional orthogonal polynomials on the given distribution.
Let the leading coefficient of $Q_{\bm{n}}(\bm{x})$ be $S_{\bm{n}}(\bm{x})$. Then for $|\bm{n}| = |\bm{n}^\prime|$
\[
\mathbb{E}\big [S_{\bm{n}}(\bm{Y})Q_{\bm{n}^\prime}(\bm{Y})\big ] = \mathbb{E}\big [Q_{\bm{n}}(\bm{Y})Q_{\bm{n}^\prime}(\bm{Y})\big ]
=\delta_{\bm{n}\bm{n}^\prime},
\]
because $Q_{\bm{n}^\prime}(\bm{Y})$ is orthogonal to polynomials in $\bm{Y}$ of degree less than $|\bm{n}^\prime|$.
Thus
\[
\mathbb{E}\big [S_{\bm{n}}(\bm{Y})Q_n(\bm{x},\bm{Y})\big ]= Q_{\bm{n}}(\bm{x}),
\]
uniquely determining $Q_{\bm{n}}(\bm{x})$ by $S_{\bm{n}}(\bm{x})$ among all orthogonal polynomial sets on $\bm{X}$. 
\qed

\bigskip

Because of Lemma 1 to check that two orthogonal polynomials sets on the same distribution are identical (up to normalizing constants) it is sufficient to check that the leading coefficients are proportional. We now calculate the leading coefficients in \cite{X2013}'s orthogonal polynomials on the multinomial and show an identity with our multivariate Krawtchouk polynomials. The notation in  \cite{X2013} is adapted to agree with notation in this paper.
%
%
In our orthogonal polynomials $Q_{\bm{n}}(\bm{x},\bm{u})$ has a single leading term proportional to
\[
S_{\bm{n}}(\bm{x})=\prod_{j=1}^{d-1}S_j^{n_j},
\text{~~~where~~~}
S_j = \sum_{k=1}^{d}u_k^{(j)}x_k,
\]
with $\{u^{(j)}\}$ a set of orthogonal functions on 
$\bm{p} = (p_1,\ldots ,p_{d})$.

We now describe the conditional binomial constructed polynomials in \cite{X2013}.
The 1-dimensional Krawtchouk polynomials there are defined by
\[
K_n(x;p,N) =~ _2F_1(-n,-x;-N;p^{-1}),\>n=0,\ldots ,N,
\]
where, with a standard definition,
\[
_2F_1(a,b;c;z) = \sum_{r=0}^\infty
\frac{a_{(k)}b_{(k)}}{c_{(k)}}\frac{z^k}{k!}.
\]
The scaling is such that $K_n(0;p,N)=1$. 
The conditional binomial multidimensional Krawtchouk polynomials are 
defined by
\begin{eqnarray}
&&K_{\bm{n}}(\bm{x};\bm{p},N)
=
\frac{ (-1)^{|\bm{n}|} }{ (-N)_{(|\bm{n}|)} }
\prod_{j=1}^{d-1}
\frac{ p_j^{n_j} }{ (1-|\bm{p}_{j-1}|)^{n_j} }
(-N + |\bm{x}_{j-1}|+|\bm{n}^{j+1}|)_{(n_j)}
\nonumber \\
&&~~~~~~~~~~~~~~~~~~~~~~
\times K_{n_j}\Big (x_j;\frac{p_j}{1-|\bm{p}_{j-1}|},N-|\bm{x}_{j-1}|-|\bm{n}^{j+1}|\Big ),
\label{cbK:0}
\end{eqnarray}
where for $j < d$, $\bm{n}^j = (n_j,\ldots,n_{d-1})$, $|\bm{n}^j| = n_{j}+\cdots + n_{d-1}$ and for notational convenience 
 $|\bm{n}^{d}| = 0$.

\bigskip

\noindent
{\bf Proposition 2.} \emph{The conditional binomial construction of multidimensional Krawtchouk polynomials, \cite{X2013}, is a special case of the multidimensional Krawtchouk polynomials where the orthogonal basis is the (unscaled) Irwin-Lancaster basis: $u^{(0)}=1$ and
\begin{equation}
u_k^{(j)} = \begin{cases}
 0&k <j,\\
 -(1-|\bm{p}_j|)/p_j&k = j,\\
 1&k = j+1,\ldots, d.
 \end{cases}
 \label{IL:0}
 \end{equation}
 for $j=1,\ldots, d-1$; $k=1,\ldots ,d$.
 }
 \bigskip
 
\noindent
{\bf Proof.} By lemma 1, it is sufficient to show that the leading terms in both sets of orthogonal polynomials are identical. Using the hypergeometric expansion
\begin{eqnarray}
&&K_{n_j}\Big (x_j;\frac{p_j}{1-|\bm{p}_{j-1}|},N-|\bm{x}_{j-1}|-|\bm{n}^{j+1}|\Big )(-N + |\bm{x}_{j-1}|+|\bm{n}^{j+1}|)_{(n_j)}
\nonumber \\
&&~~=\sum_{k=0}^{n_j}
\frac{
(-n_j)_{(k)}(-x_j)_{(k)}
}
{
(-N + |\bm{x}_{j-1}|+|\bm{n}^{j+1}|)_{(k)}
}
\Bigg (
\frac{1-|\bm{p}_{j-1}|}{p_j}
\Bigg )^k\frac{1}{k!}
\nonumber \\
&&~~~~~~~~~~~\times (-N + |\bm{x}_{j-1}|+|\bm{n}^{j+1}|)_{(n_j)}.
\label{cbK:1}
\end{eqnarray}
The leading coefficient in
\begin{eqnarray*}
&&\frac{
(-N + |\bm{x}_{j-1}|+|\bm{n}^{j+1}|)_{(n_j)}
}
{
(-N + |\bm{x}_{j-1}|+|\bm{n}^{j+1}|)_{(k)}
}
\nonumber \\
&&~~
= (-N + |\bm{x}_{j-1}|+|\bm{n}^{j+1}|+k)\cdots
(-N + |\bm{x}_{j-1}|+|\bm{n}^{j+1}|+n_j-1)
\end{eqnarray*}
is 
\[
(-N+|\bm{x}_{j-1}|)^{n_j-k}.
\]
The leading coefficient in (\ref{cbK:1}) is therefore
\begin{eqnarray*}
&&\sum_{k=0}^{n_j}
(-n_j)_{(k)}(-x_j)^k
\Bigg (
\frac{1-|\bm{p}_{j-1}|}{p_j}\Bigg )^k
(-N+|\bm{x}_{j-1}|)^{n_j-k}
\frac{1}{k!}
\nonumber \\
&&~~~
= (-1)^{n_j}\Bigg ( -x_j\frac{1-|\bm{p}_{j-1}|}{p_j} + N-|\bm{x}_{j-1}|\Bigg )^{n_j}
\nonumber \\
&&~~~
= (-1)^{n_j}\Bigg ( -x_j\frac{1-|\bm{p}_{j}|}{p_j} + x_{j+1}+\cdots + x_{d}\Bigg )^{n_j}
\nonumber \\
&&~~~
= (-1)^{n_j}\Bigg ( \sum_{k=1}^{d}u^{(j)}_kx_k\Bigg )^{n_j}.
\end{eqnarray*}
The leading coefficient in 
$K_{\bm{n}}(\bm{x};\bm{p},N)$ is then seen to be
\[
\frac{ 1 }{ (-N)_{(|\bm{n}|)} }
\prod_{j=1}^{d-1}
\frac{ p_j^{n_j} }{ (1-|\bm{p}_{j-1}|)^{n_j} }
\Bigg ( \sum_{k=1}^{d}u^{(j)}_kx_k\Bigg )^{n_j}
\]
which completes the proof.
\qed

We now calculate the proportionality constants in the two systems.
If $\bm{x} = N\bm{e}_{d}$ the generating function (\ref{s:2.4}) for the multivariate Krawtchouk polynomials is
\[
\Bigg ( 1 + \sum_{j=1}^{d-1}u_{d}^{(j)}w_j\Bigg )^{N}
=\Bigg ( 1 + \sum_{j=1}^{d-1}w_j\Bigg )^{N}
\]
so
\[
Q_{\bm{n}}\big ({N\bm{e}_{d}},\bm{u}\big )
= \frac{N!}{(N-|\bm{n}|)!\prod_{j=1}^{d-1}n_j!}.
\]
As a comparison
\begin{eqnarray*}
K_{\bm{n}}(N\bm{e}_{d};\bm{p},N) 
&=& \frac{ (-1)^{|\bm{n}|} }{ (-N)_{(|\bm{n}|)} }
\prod_{j=1}^{d-1}
\frac{ p_j^{n_j} }{ (1-|\bm{p}_{j-1}|)^{n_j} }
(-N +|\bm{n}^{j+1}|)_{(n_j)}
\nonumber \\
&=& 
(-1)^{|\bm{n}|}\prod_{j=1}^{d-1}
\Bigg (\frac{ p_j}{ 1-|\bm{p}_{j-1}| }\Bigg )^{n_j}
\end{eqnarray*}
Comparing the two polynomials at $\bm{x}=N\bm{e}_{d+1}$ gives the next corollary.
\bigskip

\noindent
{\bf Corollary 1.}
\begin{equation}
K_{\bm{n}}(\bm{x};\bm{p},N) 
= \frac{1}{(-N)_{(|\bm{n}|)}}\prod_{j=1}^{d-1}n_j!
\Bigg (\frac{ p_j}{ 1-|\bm{p}_{j-1}| }\Bigg )^{n_j}
Q_{\bm{n}}(\bm{x},\bm{u}).
\label{cbK:i}
\end{equation}
The generating function definition makes it easy to compute various transforms of $\Qn$ as a product of linear forms
\begin{equation}
\bE\Big [\prod_{i=1}^d\phi_i^{X_i}\Qn(X)\Big ] = 
{N \choose |\bm{n}|}{|\bm{n}|\choose \bm{n}}
T_0(\phi)^{N-|\bn|}
T_1(\phi)^{n_1}\cdots T_{d-1}(\phi)^{n_{d-1}}
\label{s:2.8}
\end{equation}
where
\[
T_i(\phi) = \sum_{j=1}^{d}\phi_jp_ju_{ j}^{(i)}, \>0 \leq i \leq d-1.
\]
Using random variable notation can provide elegant formulae. Let $Z_1,\ldots ,Z_N$ be independent identically distributed random variables with $P(Z=k) = p_k,\>1 \leq k \leq d$. Set $X_i = |\{k:Z_k=i\}|$ so 
$(X_1,\ldots,X_d)$ has a $m(\bx,\bp)$ distribution. Then (\ref{s:2.5}) gives (with both sides random variables)
\begin{equation}
G(\bX;\bw;\bu) = \prod_{k=1}^N\Bigg (1 + \sum_{l=1}^{d-1}w_lv_{Z_k}^{(l)}\Bigg ).
\label{s:2.9}
\end{equation}
Expanding the right hand side gives
\begin{equation}
Q_{\pmb{n}}(\pmb{X},\pmb{u}) = \sum_{\{A_l\}}\prod_{k_1\in A_1}
u_{Z_{k_1}}^{(1)}
\cdots \prod_{k_{d-1}\in A_{d-1}}
u_{Z_{k_{d-1}}}^{(d-1)},
\label{s:2.10}
\end{equation}
where the summation is over all partitions of $N$ into subsets $\{A_l\}$ such that $|A_l|=n_l$, $l=1,\ldots ,d-1$. This shows that $\Qn (\bx,\bu)$ is a polynomial of degree $n-|\bn|$ in $\big (S_1(\bx),\ldots ,S_{d-1}(\bx)\big )$ with
\[
S_i(\bx) = \sum_{j=1}^du_j^{(l)}x_j,\>i=1,\ldots ,d-1
\]
whose only term of maximal degree is $\prod_{k=1}^{d-1}S_k^{n_k}(\bx)$.

General $1$-dimensional orthogonal polynomials $\{P_n(x)\}$ satisfy a three term recurrence for $xP_n(x)$. A similar recurrence holds for $S_i(\bm{x})Q_{\bn}(\bm{x};\bu)$. Scale $Q_{\bn}(\bm{x};\bu)$ so that the leading coefficient of $\prod_{i=1}^{d-1}S_i(\bm{x})^{n_i}$ is unity by taking $Q^\ast_{\bn}(\bm{x};\bu) = \prod_{i=1}^{d-1}n_i!Q_{\bn}(\bm{x};\bu)$. Let
$c(i,l,k) = \sum_{j=1}^du_j^{(i)}u_j^{(l)}u_j^{(k)}p_j$. Then a generating function argument shows that
\begin{eqnarray}
S_i(\bm{x})Q^{\ast}_{\bm{n}}(\bm{x},\bu) &=& Q^{\ast}_{\bm{n}+\bm{e}_i}(\bm{x},\bu) + (N-|n|+1)Q^{\ast}_{\bm{n}-\bm{e}_i}(\bm{x},\bu)
\nonumber \\
&&~~~~~~~~~~~
+\sum_{l,k=1}^{d-1}c(i,l,k)Q^{\ast}_{\bm{n}-\bm{e}_l+\bm{e}_k}(\bm{x},\bu).
\label{xQn}
\end{eqnarray}
\section*{2.2 Three Examples}
This section develops three detailed examples; the first gives an `always available' basis $\{u^{(l)}\}$ for general $\bp$. This turns out to diagonize a Metropolis algorithm and satisfy a hypergroup property developed further in section 2.4. The second involves group characters and `explains' the hypergroup nomenclature; the third is a development from physics. It offers ways of generalizing the construction of section 2.1 to general space.
\bigskip

\noindent
\textbf{Example 2.1} One simple closed form choice of $\{u^{(l)}\}$ is given by Irwin-Helmert matrices \citep{L1969}.\textcolor{black}{The basis is a scaled version of (\ref{IL:0}).} Given $\bp$, let $a_i^2 = p_i$, $A_i^2 = p_d + \cdots +p_i$. Define a $d\times d$ matrix ${\cal U}$ with first row $(1,\ldots,1)$ and $i^{\text{th}}$ row
\[
\Bigg (\overbrace{0,\ldots, 0}^{d-i},\frac{ -A_{d+2-i}} { a_{d+1-i}A_{d+1-i} },
\frac{ a_{d+1-i} }{ A_{d+2-i}A_{d+1-i} },
\ldots,
\frac{a_{d+1-i}}{A_{d+2-i}A_{d+1-i}}
\Bigg ),\>2 \leq i \leq d.
\]
Thus when $d=5$, ${\cal U}$ is
\begin{equation*}
\begin{bmatrix}
1&1&1&1&1\\
0&0&0&-\frac{A_5}{a_4A_4}&\frac{a_4}{A_4A_5}\\
0&0&-\frac{A_4}{a_3A_3}&\frac{a_3}{A_3A_4}&\frac{a_3}{A_3A_4}\\
0&-\frac{A_3}{a_2A_2}&\frac{a_2}{A_2A_3}&\frac{a_2}{A_2A_3}&\frac{a_2}{A_2A_3}\\
-\frac{A_2}{a_1}&\frac{a_1}{A_1A_2}&\frac{a_1}{A_1A_2}&\frac{a_1}{A_1A_2}&\frac{a_1}{A_1A_2}
\end{bmatrix}
\end{equation*}
The rows of ${\cal U}$ are an orthonormal basis. Thus
\begin{equation}
u_j^{(0)} \equiv 1, \>
u_j^{(i)} = \begin{cases} 0&1 \leq j \leq i-1\\
-\frac{A_{i+1}}{a_iA_i}&j=i\\
\frac{a_i}{A_{i}A_{i+1}}&j > i.
\end{cases}
\label{s:2.11}
\end{equation}
These matrices were used by Lancaster and Irwin to decompose the usual chi-square test for goodness of fit to a multinomial model into $d$ orthogonal pieces. See \citet{L1969}. \citet{S2010}  has observed that the $u^{(l)}$ diagonalize a natural Markov chain. We state this formally:
\bigskip

\noindent
\textbf{Proposition 3.} (Saltzman) Let $p_1 \geq \cdots \geq p_d$ be a probability distribution on $[d] = \{1,2,\ldots ,d\}$ and $u^{(l)}$ be defined by $(\ref{s:2.11})$. Let $K(i,j)$ be the `random scan Metropolis' Markov chain on [d]: From $i$, pick $j$ uniformly in [d]. If 
$j \leq i$, move to $j$. If $j > i$, flip a $p_j/p_i$ coin. If heads, move to $j$; else stay at $i$. This is a reversible Markov chain with stationary distribution $\bp$. It has $u^{(l)}$ as (right) eigenfunctions: $Ku^{(0)} = u^{(0)}$ and $Ku^{(l)} = \beta_lu^{(l)}$ with
\[
\beta_l = 1 - \frac{A_l^2}{da_l^2},\>1 \leq l \leq d-1.
\]  
\bigskip

\noindent
\textbf{Proof} The diagonalization of the random scan Metropolis chain is a special case of a theorem of \citet{L1996}.    His result gives similar eigen-values and eigen-vectors for the Metropolis chain with general proposal.
\qed
\bigskip

The hypergroup property of an orthogonal basis $\{u^{(l)}\}$ allows delineation of all Markov chains and all (Lancaster) bivariate distributions admiting  $\{u^{(l)}\}$ as eigen-bases. Section 2.4 below shows that the multivariate Krawtchouk polynomials satisfy the hypergroup property \emph{provided} the underlying $\{u^{(l)}\}$ do. To provide examples, we now give a necessary and sufficient condition on the Irwin-Helmert basis for this property.

One way to state the property is to transform $\{u^{(l)}\}$ into an orthogonal matrix $H$ by multiplying each column by $p_j^{1/2}$; 
thus suppose $\{u^{(l)}\}$ satisfies (2.2) and define
\begin{equation}
H_{ij} = u_j^{(i-1)}p_j^{1/2} \text{~for~}i,j \in [d]. \mbox{~(Thus~$h_{1j} = p_j^{1/2}$ for $j \in [d]$)}.
\label{s:2.12}
\end{equation}
For a general orthogonal $H$, the hypergroup property is
\begin{equation}
s(j,k,l) = \sum_{i=1}^dh_{ij}h_{ik}h_{il}/h_{id}\geq 0,\> \text{for~all~} j,k,l \in [d].
\label{s:2.13}
\end{equation}
For this to be defined, $h_{id} \ne 0$ for $i \in [d]$. 
\textcolor{black}{This property is equivalent to the usual conception of a hypergroup if we take (without loss of generality) $h_{id} > 0$ because
\[
h_{ij}h_{ik} = h_{id}\sum_{l=1}^ds(j,k,l)h_{il}.
\]
That is, the product $h_{ij}h_{ik}$ can be expressed as a linear combination of terms $h_{il}$ with non-negative coefficients.}
\textcolor{black}{The property (\ref{s:2.13})} is satisfied by the Irwin-Helmert matrices and will be assumed throughout this section.
A familiar example, with $d=4$ and $p= (\frac{1}{4},\frac{1}{4},\frac{1}{4},\frac{1}{4})$ is
\begin{equation}
H = \frac{1}{2}
\begin{bmatrix} 
\phantom{-}1&\phantom{-}1&\phantom{-}1&\phantom{-}1\\
-1&\phantom{-}1&-1&\phantom{-}1\\
\phantom{-}1&\phantom{-}1&-1&-1\\
-1&\phantom{-}1&\phantom{-}1&-1\\
\end{bmatrix}
.
\label{basis:4}
\end{equation}
However the reader will find it impossible to construct a $3\times 3$ orthogonal $H$ satisfying (\ref{s:2.13}) with first row
$(\frac{1}{\sqrt{3}},\frac{1}{\sqrt{3}},\frac{1}{\sqrt{3}})$. Indeed \citet{BH2008} show that a $d\times d$ orthogonal $H$ satisfying (\ref{s:2.13}) 
with constant first row exists if and only if $d=2^k$ for some $k \geq 1$. This raises the question: what $(p_1,\ldots ,p_d)$ admit $u^{(l)}$ and 
$H$ satisfying (\ref{s:2.13})? The proposition below gives examples if $\bm{p}$  satisfies a monotonicity property.
Say that $\bp$ is \emph{strongly monotone} if 
\begin{equation}
p_d \leq p_{d-1},\> p_d + p_{d-1} \leq p_{d-2},\ldots ,p_d+\cdots + p_2 \leq p_1.
\label{s:2.14}
\end{equation}
For example when $d=3$, $(\frac{1}{2},\frac{1}{3},\frac{1}{6})$ and $(\frac{1}{2},\frac{1}{4},\frac{1}{4})$ are strongly monotone.
\bigskip

\noindent
\textbf{Proposition 4.} Suppose that $\bp$ is strongly monotone. For the Irwin-Helmet matrix (\ref{s:2.11}) define an orthogonal matrix $H$ by (\ref{s:2.12}). Then $H$ satisfies the hypergroup property. Further the strongly monotone $\bp$ form a full dimensional (dimension $d-1$) compact-convex simplex with extreme points 
\[
(0,\ldots,0,1), 
(0,\ldots,0,\frac{1}{2},\frac{1}{2}), 
(0,\ldots,0,\frac{1}{4},\frac{1}{4},\frac{1}{2}),
\ldots,
 (\frac{1}{2^{d-1}},\frac{1}{2^{d-1}},\frac{1}{2^{d-2}},\ldots ,\frac{1}{2}).
\]
\bigskip

\noindent
\textbf{Proof} It may help the reader to have an example. For $d=5$, $a_i^2 = p_i$, $A_i^2 = p_d + p_{d-1}+\cdots +p_i$, (so $A_1=1$),
\begin{equation}
H =
\begin{bmatrix}
a_1&a_2&a_3&a_4&a_5\\
0&0&0&-\frac{A_5}{A_4}&\frac{a_4a_5}{A_4A_5}\\
0&0&-\frac{A_4}{A_3}&\frac{a_3a_4}{A_3A_4}&\frac{a_3a_5}{A_3A_4}\\
0&-\frac{A_3}{A_2}&\frac{a_2a_3}{A_2A_3}&\frac{a_2a_4}{A_2A_3}&\frac{a_2a_5}{A_2A_3}\\
-A_2&\frac{a_1a_2}{A_2}&\frac{a_1a_3}{A_2}&\frac{a_1a_4}{A_2}&\frac{a_1a_5}{A_2}
\end{bmatrix}
.
\end{equation}
Observe that $s(j,k,l)$ is invariant under permutation of indices, so without loss of generality take $j \leq k \leq l$. The \textcolor{black}{following} argument shows that 
\begin{eqnarray}
&&s(j,k,l) \geq 0\text{~if~any~of~}j,k,l =d,
\nonumber \\
&&~\text{~indeed~}s(j,j,d) = s(d,d,d)=1,\> s(j,d,d) = s(j,k,d)=0;
\label{s:2.15}
\\
&&s(i,i,j) \geq 0;
\label{s:2.16}
\\
&&s(i,j,j) =0;
\label{s:2.17}
\\
&&s(j,k,l) = 0\text{~for~}1 \leq j < k < l < d.
\label{s:2.18}
\end{eqnarray}
Inequalities (\ref{s:2.15})-(\ref{s:2.18}) hold for any $a_i$. Strong monotonicity is not needed. This covers all choices except $s(i,i,i), i \in [d]$.
The positivity of these quantities is in a one to one correspondence with the linear inequalities $p_d+\cdots + p_{i+1} \leq p_i$, $1 \leq i \leq d-1$. To prove (\ref{s:2.15}) observe \emph{eg} 
$s(j,k,d) = \sum_{i=1}^dh_{ij}h_{ik} = \delta_{jk} \geq 0$. The other cases are similar. Thus without loss of generality $1 \leq j \leq k \leq l < d$ for the rest of the proof. For (\ref{s:2.16}), observe that for
$i < j$, $s(i,i,j)$ is the sum of positive terms.

The proof of (\ref{s:2.17}) is by induction on $d$. It is basic to check for $d=2,3$. For larger $d$, first consider
\[
s(1,j,j) = \frac{a_1a_j^2}{a_d} - A_2\Biggl (\frac{a_1a_j}{A_2}\Biggr )^2
\frac{A_2}{a_1a_d} = 0.
\]
Next consider the matrix $H$ with the first column and last row deleted. This is a $(d-1)\times (d-1)$ matrix of the same form with the first row divided by $A_2$. By induction
\[
\frac{p_d}{A_2} \leq \frac{p_{d-1}}{A_2}, \frac{p_d + p_{d-1}}{A_2}\leq \frac{p_{d-2}}{A_2},\ldots ,\frac{p_d+\cdots +p_j}{A_2}\leq \frac{p_2}{A_2}.
\]

This finishes the proof of the hypergroup property. The final claims are about the set of strongly monotone probabilities. Observe that 
$p_d+\cdots + p_2 \leq p_1$ is equivalent to $p_1 \geq \frac{1}{2}$. From this $p_2 \leq \frac{1}{2}$. Also $p_2+p_3 \leq \frac{1}{2}$. Along with $p_3 \leq p_2$, this gives $p_2 \leq \frac{1}{4}$. Continuing gives $p_i \leq \frac{1}{2^i}, 2 \leq i \leq d$. From this the claimed extreme points are all extreme. They can be seen to be all of the extreme points with unique representation by a similar greedy algorithm.  Given strongly monotone 
$(p_1,\dots, p_d)$, subtract off $p_d2^{d-1}(\frac{1}{2^{d-1}},\frac{1}{2^{d-1}},\ldots ,\frac{1}{2})$. The difference is positive and normalizing gives a strictly monotone probability with the first entry zero. Continuing gives $\bp$ as a linear combination of extreme points. The representation is unique because the extreme points 
\[
(1,0,0,\ldots,0),(\frac{1}{2},\frac{1}{2},0,\ldots,0),(\frac{1}{2},\frac{1}{4},\frac{1}{4},0,\ldots,0),\ldots,(\frac{1}{2},\frac{1}{4},\frac{1}{8},\ldots,\frac{1}{2^{d-1}},\frac{1}{2^{d-1}})
\]
are linearly independent.
\qed
\bigskip

Proposition 4 suggests a question about the Irwin-Lancaster bases.
$\{u^{(l)}\}$ posesses the \textcolor{black}{GKS} property, \citet{BH2008}, if for all $l,r$
\[
u^{(l)}u^{(r)} = \sum_kc^k_{lr}u^{(k)},
\]
where $c_{lr}^k \geq 0$ for $l,r,k = 0,\ldots, d-1$.
The \textcolor{black}{GKS} property insures a variety of probabilistic inequalities hold. Do the Irwin-Helmert bases (\ref{s:2.11}) satisfy \textcolor{black}{GKS} under the assumptions of Proposition 4?
The next proposition affirms that this is true.
\bigskip

\noindent
\textbf{Proposition 5.} The Irwin-Lancaster bases possess the \textcolor{black}{GKS} property if and only if $\bm{p}$ is strongly monotone.
\bigskip

\noindent
\textbf{Proof} The \textcolor{black}{GKS} property holds if and only if for $0 \leq l \leq m \leq d-1$
\begin{equation}
\sum_{j=1}^du_j^{(l)}u_j^{(m)}u_j^{(r)}p_j \geq 0.
\label{GKS:0}
\end{equation}
It is convenient to define 
\begin{eqnarray}
v_j^{(r)} &=& \frac{A_rA_{r+1}}{a_r}u_j^{(r)}
\nonumber \\
&=&
\begin{cases}
-\frac{p_{r+1}+\cdots +p_d}{p_r}&\text{~if~}j=r\\
1&\text{~if~} j > r.
\end{cases}
\label{GKS:1}
\end{eqnarray}
Then (\ref{GKS:0}) holding is equivalent to 
\[
c(l,m,r) \overset{\text{Def}}= \sum_{j=1}^dv_j^{(l)}v_j^{(m)}v_j^{(r)}p_j \geq 0.
\]
Evaluating the triple sum
\begin{equation}
c(l,m,r) = 
\begin{cases}
-\frac{p_{r+1}+\cdots +p_d}{p_r}p_r + p_{r+1}+\cdots +p_d = 0,
&\text{~if~} l \leq m < r
\\
\frac{(p_{r+1}+\cdots +p_d)^2}{p_r^2}p_r+ p_{r+1}+\cdots +p_d  \geq 0,
&\text{~if~} l < m = r
\\
\frac{(p_{r+1}+\cdots +p_d)}{p_r^2}
\big (p_r^2 - (p_{r+1}+\cdots +p_d)^2\big ) \geq 0,
&\text{~if~}l=m=r.
\end{cases}
\label{GKS:2}
\end{equation}
Positivity in the last case holds if and only if $\bm{p}$ is strongly monotone.
\qed
\bigskip

\noindent
\textbf{Hypergroups for Groups}. Let ${\cal G}$ be a finite group with conjugacy classes $\C_1,\C_2, \ldots ,\C_d$ and irreducible characters
$\chi_1,\ldots ,\chi_d$. Label those so that $\C_d = \{\text{id}\}$, $\chi_1 \equiv 1$ (the trivial character). Background can be found in
\citet{JL2001},
\citet{I1994},
\citet{D1988}.
Define $p_i = |\C_i|/|{\cal G}|$, for $i=1,\ldots,d$. It is a classical fact   that 
\[
\sum_i \chi_i(\C_j)\chi_i(\C_k)\chi_i(\C_l) \geq 0\text{~for~all~}j,k,l \in [d],
\]
see \emph{eg} \citet{BH2008}, Proposition 2.6.
This implies
\bigskip

\noindent
\textbf{Proposition 6.} Let ${\cal G}$ be a finite group with $d$ conjugacy classes. Suppose that all of the characters are real valued. Then the matrix 
$H$
\[
h_{ij} = \chi_i(\C_j) \sqrt{p_j},\> (p_j = |\C_j|/|{\cal G}|)
\]
is orthogonal and satisfies the hypergroup and \textcolor{black}{GKS} properties.
\bigskip

\noindent
\textbf{Example} Let ${\cal G} = {\cal C}_2^n$, the group of binary $n$-tuples under commutative addition. This is an abelian group so the conjugacy classes are single points. For $x\in \G$, let $\chi_x(y) = (-1)^{x\cdot y}$. $\{\chi_x\}_{x \in \G}$ are the irreducible characters and $\G$ is real.
Here $|\G|=2^n$, and $p_i = 1/2^n$. When $n=2$ the relevant basis is displayed in (\ref{basis:4}).
\bigskip

\noindent
\textbf{Example} Let ${\cal G} = S_3$, the symmetric group.
There are three conjugacy classes
\[
\C_1 = \{(1,2),(1,3),(2,3)\},\> \C_2 = \{(1,2,3),(1,3,2)\},\> \C_3 = \{\text{id}\}.
\]
Thus $p_1 = \frac{1}{2}$, $p_2 = \frac{1}{3}$, $p_3 = \frac{1}{6}$. The character table and associated $H$ are (bordered rows and columns)
\begin{equation}
\kbordermatrix{
&\phantom{-}\C_1&\phantom{-}\C_2&\phantom{-}\C_3\\
\chi_1&\phantom{-}1&\phantom{-}1&\phantom{-}1\\
\chi_2&\phantom{-}0&-1&\phantom{-}2\\
\chi_3&-1&\phantom{-}1&\phantom{-}1
}
\>\>\>\>\>\>
H =
\begin{bmatrix}
\phantom{-}\frac{1}{\sqrt{2}}&\phantom{-}\frac{1}{\sqrt{3}}&\phantom{-}\frac{1}{\sqrt{6}}\\
\phantom{-}0&-\frac{1}{\sqrt{3}}&\phantom{-}\frac{2}{\sqrt{6}}\\
-\frac{1}{\sqrt{2}}&\phantom{-}\frac{1}{\sqrt{3}}&\phantom{-}\frac{1}{\sqrt{6}}
\end{bmatrix}
.
\end{equation}
\bigskip

The orthogonal matrices $H$ in these examples also have a probabilistic interpretation; they are the eigenvectors of any random walk on $\G$ which is constant on conjugacy classes, see \citet{D1988}, chapter 3.
\bigskip

\noindent
\textbf{Example} Bosonic Fock space and second quantization.

The multivariate Krawtchouk polynomials are closely related to a basic construction in modern physics. This connection illuminates the construction, suggesting natural generalizations to infinite spaces and to varying numbers of particles. Good references for physics are \citet{D2010}, chapter 3, \citet{RS1975}, chapter X.7 and \citet{F1972}, chapter 6.7.

We begin with the general story -- symmetrized tensors, then symmetrized to the multinomial, finally discussions and generalizations.

Let $V$ be a vector space $V^{\otimes N}$ the $N$-fold tensor product and $V_S^{\otimes N}$ the elements in
 $V^{\otimes N}$ invariant under the symmetric group $S_N$.
If $<\cdot | \cdot>$ is an inner product on $V$ then $V^{\otimes N}$ becomes an inner product space with $<a_1\otimes\cdots \otimes a_N\mid b_1\otimes\cdots \otimes b_N>=<a_1|b_1>\cdots <a_N|b_N>$. Further, $V_S^{\otimes N}$ inherits an inner product. If $\overline{a_1\otimes\cdots \otimes a_N}$ denotes symmetrization,
\[
<\overline{a_1\otimes\cdots \otimes a_N}\mid \overline{b_1\otimes\cdots \otimes b_N}> = 
\text{per~} 
\begin{bmatrix}
<a_1\mid b_1>\ldots, <a_1\mid b_N>\\
\cdots\\
<a_N\mid b_1>\ldots, <a_N\mid b_N>
\end{bmatrix},
\]
\cite{F1972}, (6.3). If $\{u^{(l)}\}_{l\in {\cal L}}$ is an orthogonal basis for $V$ (dim $V = \infty$ is allowed,) then $V^{\otimes N}$ has $u^{(l_{1})}\otimes\cdots \otimes u^{(l_{ N})}=u^{(\bm{l})}$ as an orthogonal basis with $l_i \in {\cal L}$. Symmetrizing these gives
\[
\overline{u}^{(\bm{l})} = \sum_{\sigma\in S_N}u^{(l_{\sigma 1})}\otimes\cdots \otimes u^{(l_{\sigma N})}.
\]
Because of the symmetry, $\bar{u}^{(\bm{l})}$ only depends on $n_{i_1},n_{i_2},\ldots$, with $n_l$ the number of $i$ such that $l_i = l$. The $\bar{u}^{(\bm{l})}$ are an orthonormal basis for $V_S^{\otimes w}$. In the physics literature they are often denoted by $|n_1,n_2 \cdots>$. 
Specialize to the case where $V=L^2(\mu)$ for $\mu$ a probability measure on $(\mathfrak{X},{\cal B})$. Let $\{u^{(\bm{l})}\}$ be an orthogonal basis with $u^{(0)}(x)\equiv 1$. The basis elements are 
\[
\bar{f}^{(\bm{l})}(y_1,y_2,\ldots ,y_N)
=\sum_{\sigma \in S_N}\prod_{i=1}^Nf^{(l_i)}(y_{\sigma(i)}).
\]
The degree of $f^{(\bm{n})}$ is $|\bm{n}| = n_1 + n_2 + \cdots$. Thus the degree $0$ element is 
\[
\bar{f}^{(0)} \equiv 1.
\]
The degree one basis vectors are (up to the proportionality constant $(N-1)!$)
\[
\bar{f}^{(\bm{l})}(y_1,\dots,y_N) = \sum_{\sigma \in S_N}
f^{(l_i)}(y_{\sigma(i)}).
\]
%
%

Now suppose that the underlying space $\mathfrak{X}=[d]$ with $\mu = \bm{P}$. If $\{u^{(l)}\}$ is chosen as in (\ref{s:2.2}) and $z_1,\ldots ,z_N$ has $i$ appearing $x_i(\bm{z})$ times
\[
\bar{f}^{(e_l)} = \sum_{j=1}^du_j^{(l)}x_j,\>\>1 \leq l \leq d-1.
\]
These are the linear Krawtchouk polynomials. Similarly, the higher degree basis terms are the multivariate Krawtchouk polynomials.

In the physics literature if $A$ is a self adjoint operator on $V$, then
\[
(A\otimes I\otimes\cdots \otimes I)+ (I\otimes A\otimes \cdots \otimes) + \cdots + (I\otimes I\otimes \cdots \otimes A)
\]
operates on $V_S^{\otimes N}$. It is called the second quantization of $A$. If $V=L^2(\mu)$ is our $d$-dimensional space and $A$ is the transition matrix of a $p$-reversible Markov chain, the second quantization (divided by $N$) is just ``pick a co-ordinate at random, if it is colour $i$ change it to $j$ with probability $A(i,j)$''. See \citet{F1972}, section 6.8, for the physics version. 

The development above shows how to generalize from $[d]$ to a general space. The physics development has an additional feature; the creation and destruction operators $a^+(\cdot)$ and $a^-(\cdot)$. These translate into ``add or subtract a ball'' in the multinomial picture. They necessitate working in the enhanced state space $\bigoplus_{N=0}^\infty V_S^{\otimes N}$ (Bosonic Fock space). We will not develop the story further, but believe there is a lot to be done translating between fields.
\section*{2.3 Self Duality} The univariate Krawtchouk polynomials $Q_n$ satisfy the useful duality equation $Q_n(x) = Q_x(n)$. This section defines a $d$-dimensional extension, used in section 2.4 to prove the hypergroup property which is a crucial ingredient in section three.

Duality is the easiest to describe by considering a more general class of polynomials $\widehat{Q}_{\bn^+}(\bx,H)$, where $H$ is a $d\times d$ orthogonal matrix, $\bn^+ = (n_1^+,\ldots,n_d^+)$ with $|\bn^+| = N$, $\bx = (x_1,\ldots ,x_d)$ has $|\bx| = N$.
Define
 $\wh{Q}_{\pmb{n}^+}(\pmb{x},H)$ as the coefficient of 
\begin{equation}
{N\phantom{+}\choose \pmb{n}^+}
\prod_{j=1}^dw_j^{n^+_j}
{N\choose \pmb{x}}\prod_{i=1}^dz_i^{x_i}
\label{s:2.19}
\end{equation}
in
\begin{equation}
\Big [\sum_{i,j=1}^dh_{ij}w_iz_j\Big ]^N.
\label{double:gf}
\end{equation}
 $\wh{Q}_{\pmb{n}^+}(\pmb{x},H)$ is also the coefficient of  
\[
{N\phantom{+}\choose \pmb{n}^+}
\prod_{j=1}^dw_j^{n^+_j}\text{~~~in~~~}
\prod_{j=1}^d\Bigg \{\sum_{i=1}^dh_{ij}w_i\Bigg \}^{x_i}.
\]
or the coefficient of 
\[
\phantom{+}\phantom{+}{N\choose \pmb{x}}\prod_{i=1}^dz_i^{x_j}
\text{~~~in~~~}
\prod_{i=1}^d\Bigg \{\sum_{j=1}^dh_{ij}z_j\Bigg \}^{n^+_i}.
\]
There is an evident duality
\begin{equation}
\wh{Q}_{\pmb{n}^+}(\pmb{x},H)=
\wh{Q}_{\pmb{x}}(\pmb{n}^+,H^T),
\label{duality:0}
\end{equation}
where $H^T$ denotes the transpose of $H$.
\textcolor{black}{The generating function (\ref{double:gf}) is a generating function for both systems in (\ref{duality:0}). The variable is $\pmb{x}$ with index $\pmb{n}^+$ for the system on the left side, and the variable is 
$\pmb{n}^+$ with index $\pmb{x}$ for the dual system on the right side.}
To make the connection to Krawtchouk polynomials for 
$\{u^{(l)}\}$ set
\begin{equation}
h_{ij} = u_j^{(i-1)}p_j^{1/2},\>i,j\in \{1,\ldots ,d\},
\>\>\text{~so~}h_{ij}=p_j^{1/2},
\label{s:2.20}
\end{equation}
and 
\begin{equation*}
\pmb{n}^+ = (N-|\bm{n}|,n_1,\ldots, n_{d-1}),\> |\pmb{n}| \leq N.
\end{equation*}
Then
\begin{equation}
\wh{Q}_{\pmb{n}^+}(\pmb{x},H)={N\phantom{+}\choose \pmb{n}^+}^{-1}Q_{\pmb{n}}(\pmb{x},\pmb{v})\prod_{j=1}^dp_j^{x_j/2}.
\label{s:2.21}
\end{equation}
with $H$ from (\ref{s:2.20}). $\{\wh{Q}_{\pmb{n}^+}\}$ satisfy the orthogonality relations
\begin{eqnarray}
\sum_{\pmb{x}} \wh{Q}_{\pmb{m}^+}(\pmb{x},H)\wh{Q}_{\pmb{n}^+}(\pmb{x},H){N\choose \pmb{x}} &=& \delta_{\pmb{m}\pmb{n}}{N\phantom{+}\choose \pmb{n}^+}^{-1}
\nonumber \\
\sum_{\pmb{n}^+} \wh{Q}_{\pmb{n}^+}(\pmb{x},H^T)\wh{Q}_{\pmb{n}^+}(\pmb{y},H^T){N\phantom{+}\choose \pmb{n}^+} &=& \delta_{\pmb{x}\pmb{y}}{N\choose \pmb{x}}^{-1}.
\label{s:2.22}
\end{eqnarray}
This shows that $\{\wh{Q}_{\pmb{n}^+}\}$ are orthogonal polynomials for the flat multinomial $\{p_i = 1/d\}$.
A symmetrized product form is
\begin{equation*}
\wh{Q}_{\pmb{n}^+}(\pmb{x},H) = {N\phantom{+}\choose \pmb{n}^+}^{-1}\sum_{\{A_l\}_{l=1}^d}\prod_{k\in A_1}
h_{1z_k}
\cdots \prod_{k\in A_{d}}
h_{dz_k},
\end{equation*}
where the summation is over all partitions of subsets of $\{1,\ldots ,N\}$ ,
$\{A_l\}$ such that $|A_l| = n_l$, $l = 1,\ldots ,d$ with $\{z_1,\ldots,z_N\}$ a multi-set containing $x_j$ entries equal to $j$, $j=1,\ldots ,d$. 

\section*{2.4 The Hypergroup Property}
Univariate \KP   $~Q_n(x)$ with 
$\bE\big [Q_m(X)Q_n(X)\big ] = \delta_{mn}h_n$ satisfy the hypergroup property
\[
\sum_{n=0}^Nh_nQ_n(x)Q_n(y)Q_n(z) \geq 0,\>\text{~for~all~}
x,y,z=0,1,\ldots ,N.
\]
This property was discovered and exploited by \citet{E1969} in his solution of the Lancaster problem for the Binomial distribution. \textcolor{black}{\citet{V1971} and \citet{DR1974} study group theoretic properties of the hypergroup property of the Krawtchouk polynomials.} In \citet{DG2012} the hypergroup property is used to characterize reversible Markov chains with (univariate) Krawtchouk polynomials as eigen-functions. Following Eagleson's work, a host of univariate orthogonal polynomials have been shown to satisfy the hypergroup property. A wonderful survey of this work is given by \citet{BH2008}. The purpose of this subsection is to study the hypergroup property for multivariate Krawtchouk polynomials. Applications are in section
 three. 
 The main result shows that the hypergroup property is equivalent to a hypergroup property for the chosen underlying basis $\{u^{(l)}\}$. This allows the examples developed in section 2.2 to be used for the multinomial.

It is convenient to consider scaled \MKP
\begin{equation}
Q^\diamond_{\pmb{n}}(\pmb{x},\pmb{u}) = \frac{Q_{\pmb{n}}(\pmb{x},\pmb{u})}{Q_{\pmb{n}}(N\pmb{e}_d,\pmb{u})}.
\label{s:2.23}
\end{equation}
Recall that $\bm{u} = \{u^{(0)},\ldots ,u^{(d-1)}\}$ is an orthonormal basis for functions on $[d]$ with respect to $\bm{p}$. The scaling gives 
$Q^\diamond_{\pmb{n}}(N\pmb{e}_d;\pmb{u})=1$ (the choice of coordinate $d$ in $e_d$ is chosen without loss of generality).
From the generating function (\ref{s:2.4}) 
\begin{equation}
\Qn(Ne_d,\bm{u}) = {N\phantom{+}\choose \bn^+}\prod_{i=1}^{d-1}b_i^{n_i},
\label{2.24}
\end{equation}
where
\begin{equation}
b_i = u^{(i)}_d \ne 0
\label{s:2.25}
\end{equation}
is an assumption in force throughout. This is automatic if Irwin-Helmert matrices are used. In general
\begin{equation}
\mathbb{E}\Bigl [Q^\diamond_{\pmb{m}}(\pmb{X},\pmb{u})Q^\diamond_{\pmb{n}}(\pmb{X},\pmb{u})\Bigr ] = \delta_{\pmb{m}\pmb{n}}{N\phantom{+} \choose \pmb{n}^+}^{-1}\prod_{i=1}^{d-1}b_i^{-2n_i}.
\label{s:2.26}
\end{equation}
Denote
\begin{equation}
h_{\pmb{n}}^\diamond = {N\phantom{+}\choose \pmb{n}^+}\prod_{i=1}^{d-1}b_i^{2n_i}.
\label{s:2.27}
\end{equation}
The appropriate hypergroup property is
\begin{equation}
\sum_{\pmb{n}}Q_{\pmb{n}}^\diamond(\pmb{x},\pmb{u})Q_{\pmb{n}}^\diamond(\pmb{y},\pmb{u})Q_{\pmb{n}}^\diamond(\pmb{z},\pmb{u})h_{\pmb{n}}^\diamond \geq 0.
\label{s:2.28}
\end{equation}
The first result shows that (\ref{s:2.28}) is equivalent to a similar property for the original basis $\{u^{(l)}\}$. Define an orthogonal matrix as in (\ref{s:2.12}) by
\begin{equation}
h_{ij} = u^{(i-1)}_j\sqrt{p_j},\>i,j \in [d]
\label{s:2.29}
\end{equation}
and set 
\begin{equation}
\mathfrak{s}(j,k,l) = \sum_{i=1}^dh_{ij}h_{ik}h_{il}h_{id}^{-1}.
\end{equation}
\bigskip

\noindent
\textbf{Proposition 7.} For an orthonormal basis $\{u^{(l)}_j\}$ as in (\ref{s:2.2}) with $u^{(l)}_d \ne 0$ for $l=0,1,\dots, d-1$, the hypergroup property (\ref{s:2.28}) holds if and only if 
\begin{equation}
\mathfrak{s}(j,k,l) \geq 0\text{~for~all~}j,k,l \in [d].
\label{s:2.31}
\end{equation}
\bigskip

\noindent
\textbf{Proof}
A generating function proof now follows. Note that $\wh{Q}_{\pmb{n}^+}(\pmb{x},H)=\wh{Q}_{\pmb{x}}(\pmb{n}^+,H^T)$ is the coefficient of 
\[{N\choose \pmb{x}}\prod_{j=1}^dz_j^{x_j}\]
in
\[
\prod_{i=1}^d\Bigl (\sum_{j=1}^dh_{ij}z_j\Bigr )^{n_i^+}.
\]
Multiply the sum of the triple products in  (\ref{s:2.28}) by
\begin{equation}
{N\choose \pmb{x}}\prod_{i=1}^d\alpha_i^{x_i}
{N\choose \pmb{y}}\prod_{j=1}^d\beta_i^{y_j}
{N\choose \pmb{z}}\prod_{j=1}^d\gamma_i^{z_k}
\label{dummy:0}
\end{equation}
and sum over $\pmb{x},\pmb{y},\pmb{z}$ to obtain
\begin{equation}
\Bigl [\sum_{j,k,l=1}^ds(j,k,l)\alpha_j\beta_k\gamma_l\Bigr ]^N.
\label{genpf}
\end{equation}
The coefficients of (\ref{dummy:0})
in (\ref{genpf}) are non-negative if and only if (\ref{s:2.28}) holds. The sufficiency clearly holds. For the necessity first note that for $k ,j \ne d$
\[
\mathfrak{s}(d,k,k) = 1,\> \mathfrak{s}(d,d,k)=1,\> \mathfrak{s}(d,k,j) = 0,\> \mathfrak{s}(d,d,d) = 1.
\]
The coefficient of $\alpha_j\beta_k\gamma_l(\alpha_d\beta_d\gamma_d)^{N-1}$ in (\ref{genpf}) is 
$
N\mathfrak{s}(j,k,l),
$
which is necessarily non-negative.
\qed

\section*{2.5 A Linearization Formula}
Linearization formulas for classical orthogonal polynomials express the product of two polynomials as a linear combination 
\[
P_iP_j = \sum_{k=0}^{i+j}L_{ijk}P_k.
\]
Positivity and integrality of the $L_{ijk}$ is of particular interest.
Background, motivation and references are in \citet{I2005} or \citet{HBR2000}. The celebrated Littlewood-Richardson rule \citep{M1998} gives a multivariate example. It expresses the product of two Schur functions. The hypergroup property allows such a result for multivariate Krawtchouk polynomials.

With notation as in section 2.3, let $\Qn^\diamond$ be defined by
(\ref{s:2.23}). Set
\[
\phi_{\bx\by} = m(\bm{z};\bp)\sum_{\bn} \Qn^\diamond(\bx;\bu) 
\Qn^\diamond(\by;\bu) 
\Qn^\diamond(\bm{z};\bu)h_{\bn}^\diamond.
\]
From Proposition 7, $\phi_{\bx\by}(\bm{z})$ is a probability distribution in $\bm{z}$ if and only if
 $\mathfrak{s}(j,k,l) \geq 0$ 
 for all $j,k,l\in [d]$. If this is true (\emph{ie} when the hypergroup property holds for the original basis $\{u^{(l)}\}$) then
\begin{equation}
\Qn^\diamond(\bx,\bu)\Qn^\diamond(\by,\bu) = \mathbb{E}_{\phi_{{\bx}{\by}}}\Bigl [\Qn^\diamond(\bm{Z},\bm{u})\Bigr ].
\label{s:2.26a}
\end{equation}
Expanding the right hand side gives a positive linearization formula.
\section*{3. Markov Chains with Multivariate Krawtchouk Polynomial Eigen-vectors}
This section gives many natural examples of Markov chains with multivariate Krawtchouk polynomial eigenvectors. Section 3.1 reviews the work of \citet{KZ2009} and \citet{ZL2009} on composition Markov chains. These include generalized Ehrenfest urns, chains occurring in the evolution of DNA, neutral theory of biodiversity and others. These authors have used the polynomials to get sharp rates of convergence and to build martingales to calculate moments of coalescent times.

Section 3.2 develops a non-reversible theory using bi-orthogonal expansions. This is applied to generalizations of an urn model of \citep{M2010,M2011}. Section 3.3 offers further generalizations all of which are diagonalized by multivariate Krawtchouk polynomials.  Section 3.4 applies this construction to give a universal property of the Krawtchouk construction (and many more examples).
\section*{3.1 Composition Markov Chains}
Begin with a Markov chain $\bm{P}(i,j)$ with stationary distribution 
$\bp$ on $[d]$. This induces a variety of Markov chains on the product space $[d]^N$. One may `pick a coordinate at random and change that coordinate from $i$ to $j$ with probability $\bm{P}(i,j)$'. One may change all of the coordinates in turn, independently with $\bm{P}$. More generally, one may pick a subset $S \subseteq [N]$ with probability $\mu(S)$ which is exchangeable (so $\mu(S)$ only depends on $|S|$) and change the values of coordinates in $S$ independently with $\bm{P}$.
The symmetric group $S_N$ acts on $[d]^N$ by permuting coordinates. The orbit of a point 
$\bm{x} \in [d]^N$ is $\bn = (n_1,\ldots ,n_d)$ with 
$\bm{n}_j (\bm{x}) = \#\{i:\bx_i = j\}$. 
If $\mathbb{P}_\mu(\bx,\by)$ is the chain constructed above on $[d]$ then 
$\mathbb{P}_\mu(\bx,\by) = \mathbb{P}_\mu(\bx^\sigma,\by^\sigma)$
for all $\sigma \in S_N$. It follows from Dynkins Criteria \citep{KS1976,BDPX2005} that $\mathbb{P}_\mu$ induces a Markov chain on the orbit space 
\[
\chi(d,N) = \{n_1,\ldots ,n_d: 0 \leq n_i, \sum_{i=1}^dn_i = N\}
\>\text{~so~}|\chi(d,N)| = {d+N-1\choose N}
\]
with a multinomial $m(\bx,\bp)$ distribution. These are the composition chains of \citet{KZ2009}, \citet{ZL2009} and \citet{KM2013}. These authors give a detailed development of many examples using the \MKP to give sharp rates of convergence to stationarity and to build martingales used to bound first hitting times.

\noindent

To whet the reader's appetite, here is a brief list.
\bigskip

\noindent
\textbf{Example (Ehrenfest Urns)} Consider $N$ labeled balls distributed in $d$ urns. At each stage, a set of $s$ balls (say $|s| = k$ is fixed) is chosen uniformly at random and, for each ball, if in urn $i$, it is moved to urn $j$ with probability $P(i,j)$. Some special cases due to Mizukawa are considered in section 3.2.
\bigskip

\noindent
\textbf{Example (Hoare-Rahmann chain)} With $N$ balls in $d$ urns, attempt to move all $N$ each time as follows; balls in urn $i$ are left fixed with probability $\alpha_i$ and moved with probability $1 - \alpha_i$; if moved, they are re-distributed with probability $(\theta_1,\ldots, \theta_d)$. 
\bigskip

\noindent
\textbf{Example (Evolution of DNA Chromosomes)}. Here a string of $N$ nucleotides labeled $\{A,T,C,G\}$ undergoes independent mutation from a fixed $4\times 4$ transition matrix (so $d=4$).
\bigskip

\noindent
\textbf{Example (Lightbulb problem)} There are $N$ light bulbs. At time $t$, choose a set $S_t$, and for $i\in S_t$, if the bulb is off, switch it on, if on, switch it off. Of interest are the total number of bulbs on at time $t$, and the first time all bulbs are off. The 
\KP are used to build martingales that give exact formulas for these distributions.
\bigskip

\noindent
\textbf{Example (Coalescence times for a multiperson random walk on a graph)} Fix a connected simple graph $\G$ with $d$ vertices,
distribute $N$ chips on the vertices.  At each time, a randomly chosen chip picks a nearest neighbour at random and moves. Of interest is the first time all the chips are at a common vertex.
\section*{3.2 Non-reversible Chains and Biorthogonal Expansions}
Markov chains with a Lancaster expansion are usually thought of as reversible. This section treats non-reversible chains using Biorthogonal expansions.
Some examples of \citep{M2010,M2011} are treated.

Let $P = (P_{ij})$ be a Markov transition matrix on $[d]$ with stationary distribution $\{p_i\}$. Suppose that $P$ is diagonalizable with left eigenvectors $\{p_i\beta_i^{(k)}\}_{k=0}^{d-1}$, right eigenvectors $\{\alpha_j^{(k)}\}_{k=0}^{d-1}$ and eigenvalues $\{\rho_k\}_{k=0}^{d-1}$ with $\rho_0 = 1$. Then $P$ has spectral representation
\begin{equation}
p_{ij} = p_j\big \{1 + \sum_{k=1}^{d-1}\rho_k\alpha_i^{(k)}\beta_j^{(k)}\big \},\>i,j=1,\ldots ,d.
\label{s:3.1}
\end{equation}
We may also call (\ref{s:3.1}) a Lancaster expansion. Examples are given below and a host of further examples of non-reversible chains with explicit real left and right eigen-vectors are in \citet{DPR2011}. Of course, without  reversibility, the eigen-vectors may take complex values and need not be othogonal, but they satisfy the biorthogonality relationship
\begin{equation}
\sum_{i=1}^dp_i\alpha_i^{(k)}\beta_i^{(l)} = \delta_{kl}.
\label{s:3.2}
\end{equation}
Define two sets of \MKP
\begin{equation}
\big \{\Qn(\bx,\bm{\alpha})\big \},\>\big \{\Qn(\bx,\bm{\beta})\big \}
\label{s:3.3}
\end{equation}
using the generating function (\ref{s:2.4}).

Any of the schemes of section 3.1 can now be used to get a Markov chain on $\mathfrak{X}(d,N)$ with $N$ particles and a multinomial stationary distribution.
 These chains have spectral expansions with respect to $\Qn(\bx,\bm{\alpha}),\Qn(\bx,\bm{\beta})$.
 For example, if each of the $N$ particles independently makes a transition at each stage, the transition matrix from $\bx \to \by$ is
\begin{equation}
m(\by;\bp)\big \{1 + \sum_{\{\bn: 1 \leq |\bn| \leq N\}}\gamma_{\bn}\Qn(\bx,\bm{\alpha})\Qn(\bx,\bm{\beta})\big \},
\label{s:3.6}
\end{equation}

where
\[
\gamma_{\bn} = \prod_{j=1}^{d-1}\rho_j^{n_j}.
\]
\bigskip

\noindent
\textbf{Example} \citet{M2010,M2011} considered $N$ balls distributed in $d$ urns arranged around a circle and moved to another urn in one of three schemes:\\

\noindent
(a) a randomly chosen different urn;

\noindent
(b) the next urn right (mod $d$); and

\noindent
(c) one of the two adjacent urns (mod $d$) with equal probability.
\bigskip

All of these examples have transition matrices for a single ball change which are circulants \citep{D1979}. Let $\bm{P}$ be a general $d\times d$ circulant transition matrix with first row $\{q_j\}_{j=1}^d$ and other rows rotated successively from the first, so the $i^{\text{th}}$ row is $\{q_{j-i}\}^d_{j=1}$ with subscripts taken mod $d$.
$\bm{P}$ is doubly stochastic, with a uniform stationary distribution on $[d]$. An eigenfunction expansion of $\bm{P}$ is
\begin{equation}
p_{ab} = \frac{1}{d}\sum_{k=0}^{d-1}\eta_k e^{2\pi ik(a-1)/d}
 e^{-2\pi ik(b-1)/d},\>\>a,b=1,\ldots,d,
 \label{s:3.7}
\end{equation}
where
\begin{equation}
\eta_k = \sum_{r=0}^{d-1}q_{r+1}e^{2\pi i rk/d}.
\label{s:3.8}
\end{equation}
Here $\bm{P}$ is reversible if and only if it is symmetric. (\emph{eg} cases (a), (c)).

Construct the \MKP (\ref{s:3.3}) by taking 
\[
\alpha^{(l)}_a=e^{2\pi il(a-1)/d},\>\>\beta^{(l)}_b=e^{-2\pi il(b-1)/d}, \>l=0,\ldots ,d-1,\>a,b \in [d].
\]
Then the representation (\ref{s:3.6}) holds with
\begin{equation}
\gamma_{\bm{n}} = \sum_{l=0}^{d-1}\eta_l\frac{n_l}{N}.
\label{s:3.9}
\end{equation}

In this circulant case, the eigenfunctions $\Qn$ are the monomial symmetric functions as recognized by Mizukawa. To see this, 
label the balls so that $Z_l = k$ if ball $j$ is in urn $k$. Then from (\ref{s:2.16}) 
\[
\Qn(\bX,\bu) = \sum_{\{A_l\}}\prod_{k_1\in A_1}e^{(2\pi i/d)(Z_{k_1}-1)}\cdots \prod_{k_{d-1}\in A_{d-1}}e^{(2\pi i(d-1)/d)(Z_{k_{d-1}}-1)}
\]
where the summation is over subsets of $[n]$, $\{A_l\}$ with $|A_l| = n_l,\>1\leq l \leq d$. Let $n_0=N - |\bn|$, so
$\{n_j\}_{j=0}^{d-1}$ is a composition of $N$. 
Regard $(Z_{k}-1)_{k=1}^N$ as a partition 
$0^{x_1}1^{x_2}\cdots (d-1)^{x_d}$. Let $\xi = e^{2\pi i/d}$ and 
$\Xi = \Big ((1)^{n_0}(\xi)^{n_1}\cdots (\xi^{d-1})^{n_{d-1}}\Big )$ then
\begin{equation}
\Qn(\bX,\bm{\alpha}) = m_n(\Xi).
\label{s:3.10}
\end{equation}
It is plausible that very similar results hold if $P$ is a ${\cal G}$-circulant \citep{D1990}.
\section*{3.3 Further Generalizations}
There is a natural generalization of the processes in sections 3.1, 3.2 that leads to Markov chains with Krawtchouk eigen-functions. Fix a basis $\{u^{(l)}\}$, orthonormal with respect to 
$\bp = (p_1,\ldots ,p_d)$.
 Let ${\cal L}(\bm{\beta},\bu) = \{(\beta_1,\ldots,\beta_{d-1})\}$ such that
\begin{equation}
K_{\bm{\beta}}(i,j) = p_j\big \{1 + \sum_{l=1}^{d-1}\beta_lu_i^{(l)}u_j^{(l)}\big \} = 0.
\end{equation}
This ${\cal L}$ is a non-empty, compact, convex set, the Lancaster set for $\{u^{(l)},\bp\}$ \citep{K1996}. 
Observe that $\bm{\beta}\in {\cal L}$ implies that $K_{\bm{\beta}}$ is a reversible Markov Kernel with $\bp$ a stationary distribution and $\{u^{(l)}\}$ as (right) eigen-functions. 
Since $\bm{\beta}$ are the eigen-values, $-1 \leq \beta_i \leq 1$. 
${\cal L}$ contains 1 in an open neighbourhood of zero and if $\bm{\beta} = (1,\ldots,1)$, $K_{\bm{\beta}}$ is allowable.
 Finally, ${\cal L}$ is closed under coordinate-wise product (Hadamard product) and so forms a commutative semi-group with $\bm{1}$ as identity. Further
 $
 K_{\bm{\beta}}K_{\bm{\gamma}}
 =K_{\bm{\gamma}}K_{\bm{\beta}}
= K_{\bm{\beta}\circ\bm{\gamma}}$.
Determining an exact description of ${\cal L}$ is an ongoing research area. See \citet{BH2008} and \citet{I2005} section 4.7 and the references in \citet{DKS2008}. If the $\{u^{(l)}\}$ satisfy the hypergroup property
\[
\sum_{l=0}^d
\frac{u_i^{(l)}u_j^{(l)}u_k^{(l)}}
{u_{i_0}^{(l)}}
 \geq 0
\]
for some fixed $i_0$ and all $i,j,k \in [d]$, then the extreme points of ${\cal L}$ are
\[
\bm{\beta} = \Biggl (
\frac{u_1^{(l)}}{u^{(l)}_{i_0}},\ldots,
\frac{u_d^{(l)}}{u^{(l)}_{i_0}}\Biggl ),\>0 \leq l \leq d-1.
\]
See \citet{BH2008} for a proof and section 2.2 for examples.

With ${\cal L}$ in focus, a general class of Markov chains with multinomial $(N,\bp)$ stationary distributions and Krawtchouk polynomial eigen-functions can be defined.

Let $\mu$ be an exchangeable probability on ${\cal L}^N$.
Define $K_\mu(\bx,\by)$ on $\mathfrak{X}(N,d)$ as the orbit chain on $[d]^N$ derived from picking $(\bm{\beta}_1,\ldots, \bm{\beta}_d)$ from $\mu$ and moving the $j^\text{th}$ coordinate with $K_{\bm{\beta}_j}(j,j^\prime)$.
The eigen-values of the product chain are
\[
\lambda_{i_1\ldots i_N} = \int \bm{\beta}_{i_1}\cdots \bm{\beta}_{i_d}
\>\mu(d{\bm{\beta}_{1}}\cdots d\bm{\beta}_{N}).
\]
Because of exchangeability, this only depends on
$\bm{i} = (i_1,i_2,\ldots, i_N)$ through
$n_1(\bm{i}), \ldots , n_d(\bm{i})$.
This defines
\begin{equation}
\lambda_{n_1\ldots n_d}.
\end{equation}
Summarizing:
\bigskip

\noindent
\textbf{Proposition 8.} Let $\{u^{(l)}\}_{l=0}^{d-1}$ be an orthonormal set on $[d]$ with repect to $\bm{p}$. Let $\mu$ be an exchangeable probability on ${\cal L}^N$. The Markov chain $K_\mu$ on
$\mathfrak{X}(N,d)$ has a $m(\bx,\bp)$ stationary distribution, multinomial Krawtchouk eigen-functions $\big \{\Qn(\bx,\bu)\big \}_{|\bn|\leq N}$ with eigenvalues $\{\lambda_{\bn}\}_{|\bn| \leq N}$.

This construction includes all of the examples in section 3.1. It shows that the \MKP have a kind of universal quality, diagonalizing the orbit chains of arbitrary products of Markov chains with $\{u^{(l)}\}$ as eigen-vectors.

The set of all exchangeable probabilities $\mu$ on ${\cal L}^N$ is a convex simplex whose extreme points are straightforward to describe: put $N$ balls in an urn labeled with $\bm{\beta}^1,\bm{\beta}^2, \ldots ,\bm{\beta}^N$ 
with $\bm{\beta}^i\in {\cal L}$ and draw them out sampling without replacement \citep{DF1980}.

\section*{Acknowledgment} Thanks to Yuan Xu for his substantial comments on the paper and to Pratha Dharmawansa for his carefull reading and list of corrections. Robert Griffiths was supported by the Department of Statistics, Stanford University in 2011; the Miller Foundation, Berkeley, in a visit to the Department of Statistics, Berkeley in 2012; the Clay Mathematics Institute in a visit to the University of Montreal in 2013; and the Institute of Statistical Mathematics, Tokyo, in 2014. He thanks the institutions for their support and hospitality.

\bigskip

\end{document}